\documentclass{gtpart}     % Basic GT/GTM/AGT style
\gtart  
\title{Commensurations and Subgroups of\\Finite Index of Thompson's Group $F$}

\author{Jos\'e Burillo}
\givenname{Jos\'e}
\surname{Burillo}
\address{Departament de Matem\`atica Aplicada IV\\
Universitat Polit\`ecnica de Catalunya\\\newline
Escola Polit\`ecnica Superior de Castelldefels\\\newline
Avda. Del Canal Olímpic 15\\
08860 Castelldefels (Barcelona)\\Spain}
\email{burillo@ma4.upc.edu}
\urladdr{http://www-ma4.upc.edu/~burillo/}

\author{Sean Cleary}
\givenname{Sean}
\surname{Cleary}
\address{Department of Mathematics R8133\\
The City College of New York\\\newline 
Convent Ave \& 138th\\ 
New York, NY 10031\\USA}
\email{cleary@sci.ccny.cuny.edu}
\urladdr{http://www.sci.ccny.cuny.edu/~cleary/}

\author{Claas E R\"over}
\givenname{Claas E}
\surname{R\"over}
\address{Department of Mathematics\\
University of Ireland, Galway\\\newline
University Road\\
Galway\\
Ireland}
\email{claas.roever@nuigalway.ie}
\urladdr{http://www.maths.nuigalway.ie/~chew/}
\volumenumber{}
\issuenumber{}
\publicationyear{}
\papernumber{}
\startpage{}
\endpage{}
\doi{}
\MR{}
\Zbl{}
\received{}
\revised{}
\accepted{}
\published{}
\publishedonline{}
\proposed{}
\seconded{}
\corresponding{}
\editor{}
\version{}

\newtheorem{thm}{Theorem}[section]    % Standard theorem environment
\newtheorem{lemma}[thm]{Lemma}          % Lemma environment with numbering 
\newtheorem{prop}[thm]{Proposition}%
\newtheorem{corol}[thm]{Corollary}%
\makeop{Com}
\makeop{Aut}
\newcommand{\com}[1]{\ensuremath{\Com(#1)}}%gives Com(?)
\newcommand{\coplus}[1]{\ensuremath{\Com^+(#1)}}%
\newcommand{\rco}[2]{\ensuremath{\Com_{#2}(#1)}}%gives Com_{?2}(?1)
\newcommand{\iv}{^{-1}}%for inverses in groups
\newcommand{\aut}[1]{\ensuremath{\Aut(#1)}}%gives Aut(?1)

\newcommand{\pl}{\ensuremath{P}}% PL maps with lots of conditions
\newcommand{\plp}{\ensuremath{P_+}}% orientation preserving such maps
\newcommand{\ganz}{\ensuremath{\mathrm{\sf Z}\!\!\!\!\:{\sf Z}}}% integers
\newcommand{\real}{\ensuremath{{\mathrm{\sf I}}\!{\mathrm{\sf R}}}}% reals
\newcommand{\ratio}{\ensuremath{\mathrm{\sf Q}\!\!\!\!\:\!{\sf I}\,\,}}% rationals
\newcommand{\dyadic}{\ensuremath{\ganz[\frac{1}{2}]}}% the dyadic integers

\newcommand{\tr}[1]{\mbox{\ensuremath{\textrm{#1}}}}%shorthand for \textrm
\newcommand{\mc}[1]{\mbox{$\mathcal{#1}$}}%shorthand for \mathcal

\begin{document}

\begin{abstract}    % type your abstract below
We determine the abstract commensurator \com{F} of Thompson's group
$F$ and describe it in terms of piecewise linear homeomorphisms of
the real line. We show
\com{F} is not finitely generated and determine which
subgroups of finite index in $F$ are isomorphic to $F$.
We also show that the natural map from the commensurator group  to the
quasi-isometry group of $F$ is injective.
\end{abstract}

\begin{asciiabstract}    % type your abstract below
We determine the abstract commensurator Com(F) of Thompson's group
F and describe it in terms of piecewise linear homeomorphisms of
the real line. We show
Com(F) is not finitely generated and determine which
subgroups of finite index in F are isomorphic to F.
We show that the natural map from the commensurator group  to the
quasi-isometry group of F is injective.
\end{asciiabstract}

\maketitle

%%%%%%%%%%%%%%%%%%%%   Start of main body of article

\section*{Introduction}

Thompson's groups have been extensively studied since their
introduction by Thompson in the 1960s,
despite the fact that Thompson's account \cite{thom} appeared only in 1980. 
They have provided examples
of infinite finitely presented simple groups, as well as some other
interesting counterexamples in group theory (see for example, Brown and Geoghegan \cite{bg}). Cannon, Floyd and Parry \cite{cfp}
give an excellent introduction to Thompson's groups
where many of the basic results used below are proven carefully.

Automorphisms for Thompson's group $F$ were studied by Brin
 \cite{brin}, where a key theorem by McCleary and Rubin
\cite{mcclearyrubin} is used to realize each automorphisms as conjugation by
a piecewise linear map. Here, we generalize from
automorphisms to commensurations, which are isomorphisms between two
subgroups of finite index. These form a group (under a natural equivalence
relation involving passing to smaller yet still finite-index subgroups), called the commensurator
group.

We classify finite-index subgroups of $F$, and then
we extend Brin's results from automorphisms to commensurations,
again realizing every commensuration as conjugation by a piecewise
linear homeomorphism of the real line. These
maps exhibit a
particular structure, satisfying an affinity condition in the
neighborhood of $\infty$ which we use to find the algebraic
structure of the commensurator of $F$.

Commensurators have proven to be an effective tool for investigating
quasi-isometries of a group to itself,  and for effectively analyzing rigidity,
particularly of lattices.  In the case of $F$, the only
quasi-isometries of $F$ known previously were automorphisms.
This paper provides a wide array of examples of
quasi-isometries, since all commensurations are quasi-isometries,
and we prove in Section 5 that the commensurator group embeds
into the quasi-isometry group in the case of $F$. 

Our approach is 
algebraic, but we note that elements of the commensurator of $F$ can 
be represented by marked, infinite, eventually periodic, binary 
tree pair diagrams. We also note that recently Bleak and Wassink
\cite{bleakwassink} have independently described the finite-index 
subgroups of $F$, using different methods.

The paper is organized as follows. In Section \ref{sec:def} we give the
necessary definitions, and in Section \ref{sec:sgps} the first basic results for
the finite-index subgroups of $F$. In Section \ref{sec:com} the main result
about the commensurator is stated and proved, and in Section \ref{sec:struct}
 its algebraic structure is given. The proof of the embedding
of the commensurator group into the quasi-isometry group is given in
Section \ref{sec:qi}. 

\subsection*{Acknowledgements}
The authors thank Matt Brin and Dan Margalit for useful conversations and are  grateful for the hospitality of the Centre de Recerca Matem\`atica. The first author acknowledges support from MEC grant \#MTM2006-13544-C02. The second author acknowledges support from the National Science Foundation and from PSC-CUNY Research Award \#69034.

\section{Definitions}\label{sec:def}

Let \pl\ denote the group of all homeomorphisms $f$ from \real\ to
itself that

\begin{enumerate}
\item are piecewise linear with a discrete (but possibly infinite) set
of breakpoints (discontinuities of the derivative of $f$),
\item use only slopes that are integral powers of $2$,
\item have their breakpoints in the set \dyadic and
\item satisfy $f(\dyadic)\subset \dyadic$.
\end{enumerate}

It is easy to check that each element $f$ of \pl\ actually satisfies
$f(\dyadic)=\dyadic$ and that \pl\ has a subgroup of index two which contains
only the order preserving elements. We denote this subgroup by \plp.
The quotient $\pl/\plp$ is generated by the image of the
homeomorphism $\tau\co t\mapsto -t$.

Let $f\in\pl$. We call $f$ {\em integrally affine}  if
$f(t)=\varepsilon t+p$ for some integer $p$ and $\varepsilon\in\{\pm
1\}$. We say $f$ is  {\em periodically affine} if
$f(t+p)=f(t)+q$ for some non-zero $p,q\in\real$
and  {\em integrally periodically affine} if $p$ and $q$ are integers.
Note that all integrally affine maps are  integrally periodically affine
with $q=\pm p$ depending on whether $f$ is in \plp\ or not. 

When $\mc{P}$ is any of the above properties, then we call $f$ 
{\em eventually $\mc{P}$} if $f$ satisfies $\mc{P}$ for all
$t\in\real$ with $|t|>M$ for some $M>0$; here $|t|$ denotes the
absolute value of $t$. For example, $f\in\plp$ is eventually integrally affine
if there exist $l,r \in\ganz$, $M\in\real$, $M>0$, so that $f(t)=t+r$ for all
$t> M$ and $f(t)=t+l$ for all $t<- M$. Notice that $l$ and $r$ may
well be different. 

It is well-known that Thompson's group $F$ is isomorphic to the subgroup of \plp\
consisting of all eventually  integrally affine elements (see \cite{cfp}).  
It is  easy to see that the commutator subgroup $F'$
of $F$ consists of all eventually trivial elements of \plp (those where
eventually $f(t)=t$). 
This group is denoted by $BPL_2(\real)$ by Brin \cite{brin},
where $B$ stands for bounded support.

\section{Finite-index Subgroups of $F$}\label{sec:sgps}

Let $f$ be an element of $F$. Since $f$ is eventually
integrally affine, there are two integers $l,r$ and a real
number $M>0$ such that $f(t)=t+r$ for $t>M$ and $f(t)=t+l$ for
$t<-M$. The two numbers $l$ and $r$ are precisely the two components of
the image of $f$ in $\ganz\times\ganz$ under the abelianization
map. The subgroups of finite index of $F$ are in  one-to-one
correspondence with those of its abelianization  $\ganz\times\ganz$  by the following result.

\begin{prop}\label{prereq}
Let $H$ be a subgroup of $F$ of finite index. Then $H$ contains
$F'$, the commutator subgroup of $F$, and hence $H$ is normal in
$F$. Moreover, $H'=F'$.
\end{prop}
\begin{proof}
Since $F$ is finitely generated, $H$ has only finitely many
conjugates in $F$ and the intersection of all of them, $K$ say, is
normal and of finite index in $F$. We consider $K\cap F'$, which is
thus normal and of finite index in $F'$.
Hence, since $F'$ is simple and infinite, we conclude that
$K\cap F'=F'$ and $F'\subset K\subset H$.

Hence $H$ is normal in $F$. The final claim follows from the fact
that $H'$ is contained in $F'$ but also characteristic in $H$ and
hence normal in $F$, whence $F' \subset H'$.
\end{proof}

From this fact we deduce that the finite-index subgroups of $F$ are
in bijection with those of $\ganz\times\ganz$. There is a
distinguished family among these---the subgroups
$p\ganz\times q\ganz$. We denote by $[p,q]$, $p,q\in\ganz$, the
preimage in $F$ under the abelianization homomorphism of the subgroup
$p\ganz\times q\ganz$ of $\ganz\times\ganz$. Thus $F=[1,1]$ and
$F'=[0,0]$.

\section{The Commensurator Group}\label{sec:com}

As mentioned before, a {\em commensuration} of a group $G$ is an
isomorphism $\alpha\co A\rightarrow B$, where $A$ and $B$ are subgroups
of $G$ of finite index. Two commensurations $\alpha$ and $\beta$ are
equivalent if they agree on some subgroup of finite index  in $G$.
In view of this, the product $\beta\circ\alpha$ of two commensurations
$$
\alpha\co A\to B\quad\mathrm{and}\quad \beta\co C\to D
$$
is defined on $\alpha\iv(B\cap
C)$. The set of all commensurations of $G$ modulo the above equivalence
relation, together with this composition, forms a
group called the {\em commensurator of $G$} which we denote by
$\com{G}$. If $G$ is a subgroup of the group $H$, then the (relative)
commensurator of $G$ in $H$, \rco{G}{H}, consists of all elements
$h$ of $H$ for which $G\cap G^h$ has finite index in both $G$ and
$G^h$; here $G^h=h^{-1}Gh$.

The main result of this paper is the following.

\begin{thm}\label{main}
The commensurator of $F$ is isomorphic to \rco{F}{\pl}, which
consists of all eventually  integrally periodically affine elements
(of \pl).
\end{thm}

The strategy of the proof is to find a large group where $F$ is a
subgroup, and in such a way that every commensuration can be seen as
a conjugation by an element of the large group. The group \pl\ plays
this role in the case of $F$.

In order to explain this strategy, we need some definitions and one of the main
results of McCleary and Rubin \cite{mcclearyrubin}. Let $(L,<)$ be a dense linear order. By
{\em interval} we mean a nonempty open interval. A subgroup $G$ of
\aut{L}\ is {\em locally moving} if for every interval $I$ there
exists a nontrivial element $g\in G$ which acts as the identity on
$L\setminus I$. Finally,  $G$ is {\em $n$-interval-transitive} if
for every pair of sequences of intervals $I_1 <\cdots < I_n$ and
$J_1 <\cdots < J_n$ there exists $g\in G$ such that $I_k^g\cap
J_k\neq \varnothing$ for $1\le k\le n$. Below, $\overline{L}$ denotes
the Dedekind completion of $L$ which is assumed to have no
endpoints.

\begin{thm}{\rm (McCleary--Rubin \cite{mcclearyrubin})}\qua\label{rmc} 
Assume $(L_i,<)$ is a dense
linear order without endpoints and let $G_i\subset\aut{L_i}$ be
locally moving and $2$-interval transitive, $i=1,2$. Suppose that
$\alpha\co G_1\rightarrow G_2$ is an isomorphism. Then there is a
monotonic bijection $\tau \co\overline{L}_1\rightarrow \overline{L}_2$
which induces $\alpha$, that is, $g^\alpha = \tau^{-1}g\tau$ for
every $g\in G_1$; and $\tau$ is unique.
\end{thm}

Being locally moving and having 2-interval transitivity are
local properties in the sense that a group inherits these from any
of its subgroups.

\begin{proof}[Proof of \fullref{main}]
View \dyadic\ as a dense linear order and $F$ as the eventually
integrally affine elements of \plp. Let $\alpha\co A\rightarrow B$ be a
commensuration of $F$. By Proposition~\ref{prereq}, both $A$ and $B$ contain $F'$ which is (obviously) locally moving and $2$-interval transitive (see \cite[Lemma 2.1]{brin}). So Theorem~\ref{rmc} tells us that $\alpha$
is induced by conjugation with a unique element of
$\mathrm{Homeo}(\real)$. This yields an injective homomorphism
$\Psi\co\com{F}\rightarrow\mathrm{Homeo}(\real)$.

Next, we show that the image of $\Psi$ is in fact contained in \pl.
By Proposition~\ref{prereq}, each commensuration of $F$ induces an
automorphism of $F'$. In other words, the image of $\Psi$ is
contained in $N_{\mathrm{Homeo}(\real)}(F')$, the normalizer of $F'$
in $\mathrm{Homeo}(\real)$. But this normalizer is equal to \pl\ by
Theorem 1 of Brin \cite{brin}. The existence and uniqueness statements in
Theorem~\ref{rmc} now imply that $\Psi$ is an isomorphism between
\com{F}\ and \rco{F}{\pl}, which proves the first part of
Theorem~\ref{main}.

Let $\alpha\in\com{F}$ and choose positive integers $p$ and $q$ so large
that $\alpha$ is defined on the subgroup $[p,q]$, that is
$[p,q]^\alpha$, the image of $[p,q]$ under $\alpha$, is contained in
$F$. By what was said above, we can view $\alpha$ as conjugation
by an element of \pl. So for $f\in [p,q]$ we  find
$f^\alpha=\alpha^{-1}f\alpha$ to be eventually  integrally affine.
Suppose for a moment that $\alpha$ is order preserving and that
$f(t)=t+kq$ for $t\gg 0$, where $k\in\ganz$. Then
$$f^\alpha (t)=(\alpha\circ f\circ\alpha^{-1})(t)=\alpha(f(\alpha^{-1}(t)))=\alpha(\alpha^{-1}(t)+kq)=t+r$$
must hold for some $r\in\ganz$. In other words,
$\alpha^{-1}(t+r)=\alpha^{-1}(t)+s$ for some integers $r$ and $s$
and all $t\gg 0$. Since $f$ was arbitrary, we may assume that $k\neq
0$, which implies that $s\neq 0$, and hence also $r\neq 0$.
Therefore $\alpha^{-1}$, and hence $\alpha$, must be
 integrally periodically affine near infinity. A similar calculation holds for $t\ll 0$ and also when
$\alpha$ is order reversing. Consequently, each commensuration of
$F$ must be eventually integrally periodically affine.

It remains to show that each eventually  integrally periodically
affine  $\beta\in\pl$ induces a commensuration of $F$ by
conjugation. Suppose $\beta(t+p)=\beta(t)+q$ for $t\gg 0$ and
$\beta(t+p')=\beta(t)+q'$ for $t\ll 0$, with
$p,q,p',q'\in\ganz\setminus\{0\}$. Let $U=[p',p]$ if $\beta$ is
order preserving and set $U=[p,p']$ otherwise. Then for $f\in U$, we
have
$$f^\beta(t)=\left\{\begin{array}{ll}
\beta(\beta^{-1}(t)+kp)=t+kq, & t\gg 0\\
\beta(\beta^{-1}(t)+k'p')=t+k'q', & t\ll 0\end{array}\right.$$ where
$k,k'\in\ganz$ depend on $f$. Together with a similar argument for
$\beta^{-1}$ one easily sees that $U^\beta = [q',q]$ or $[q,q']$,
depending on whether $\beta$ is order preserving or not.
Theorem~\ref{main} is thus established.
\end{proof}

We immediately obtain the following  corollaries from this result.

\begin{corol}
A subgroup $U$ of $F$ of finite index is isomorphic to $F$ if and
only if $U=[p,q]$ for some positive integers $p$ and $q$.
\end{corol}

\begin{proof} Suppose $U$ is a subgroup of finite index in $F$. If $U$ is
isomorphic to $F$, then there exists an eventually  integrally
periodically affine $\alpha\in\pl$ with $F^\alpha = U$ and
calculations as above show that $U$ must be of the form $[p,q]$. On
the other hand, the final paragraph of the proof of the theorem read with $p=p'=1$ shows that
$[q',q]$ is isomorphic to $F$ for every choice of positive integers
$q$ and $q'$. This completes the proof. \end{proof}

Finally, since each subgroup of finite index in $F$ contains $[p,q]$
for some positive integers $p$ and $q$ by Proposition~\ref{prereq},
we have the following results.

\begin{corol} Every finite-index subgroup of $F$ is
virtually $F$.
\end{corol}

\begin{corol}
A group is commensurable with $F$ if and only if it is a finite
extension of $F$.
\end{corol}

\section{The Structure of {$\mathrm{Com}(F)$}}\label{sec:struct}

Descriptions of  elements of \com{F} as conjugations in \pl\ allow us to
study its structure as a group. An element $\alpha$ of \com{F} is
eventually integrally periodically affine, so there exist positive
integers $p,p',q,q'$ and a real number $M$ such that 
\begin{gather*}
\alpha(t+p)=\alpha(t)+q,\tr{ for }t>M\\
\alpha(t+p')=\alpha(t)+q',\tr{ for }t<-M.
\end{gather*}
We need a lemma about affine functions, whose proof is elementary
and left to the reader.
\begin{lemma}
Let $f\co\real\to \real$ be an integrally periodically affine
map, and assume that there are integers $i,i',j,j'$ such that for
all $t\in\real$ we have
$$f(t+i)=f(t)+j\quad\mathrm{and}\quad f(t+i')=f(t)+j'.$$
Then we have
\begin{gather*}
f(t+r)=f(t)+s,\\
\tag*{\text{where}}r=\gcd(i,i')\quad\mathrm{and}\quad s=\gcd(j,j').
\end{gather*}
Furthermore, we have
$$
\frac{i}{j}=\frac{i'}{j'}.
$$
\end{lemma}

From this lemma, we see that the integers $p,p',q,q'$ for
 element of \com{F}  depend
only on the element.

We recall that \com{F} has a subgroup of index 2, denoted \coplus{F},
formed by the commensurations induced by conjugations by
piecewise-linear maps which preserve the orientation of $\real$.

\begin{prop} There exists a surjective homomorphism $\Phi\co\coplus{F}\to\ratio^*\times\ratio^*$ defined by 
$$\Phi(f)=\left(\frac pq,\frac{p'}{q'}\right).$$
\end{prop}
Here $\ratio^*$ denotes the multiplicative group of the positive
rational numbers.

The map is obviously well-defined due to the lemma above, and it is
very easy to see that it is a homomorphism of groups. The two components
of the map capture the behavior at both ends, eventually near $-\infty$
and eventually near $+\infty$. The two numbers $p/q$ and $p'/q'$
measure the ``rate of growth" of the map at both ends.

A corollary of this result is that, as expected, $\com{F}$ is
infinitely generated.

\section{Commensurations as Quasi-isometries}\label{sec:qi}

Let $G$ be a finitely generated group.
Quasi-isometries of $G$ can be naturally
composed, and there is a natural notion of equivalence class of
quasi-isometries. Two quasi-isometries are considered
equivalent if they are a bounded distance apart in the sense
that  $f$ and $g$ are considered
equivalent if there exists a number $M>0$ such that $d(f(t),g(t))\le
M$ for all $t$ in $G$.

Equivalence classes of quasi-isometries form elements
of the group of quasi-isometries
$QI(G)$ of $G$.
It is well known that the commensurator group admits a map to the
quasi-isometry group, since all commensurations give maps
between finite index subgroups which are canonically quasi-isometric
to the ambient group. The result we want to prove in this
section is that for Thompson's group $F$, this map is one-to-one.

\begin{thm} \label{thm:injective}The natural homomorphism $\displaystyle \com{F}\to QI(F)$ is injective.
\end{thm}

We begin with an elementary lemma.

\begin{lemma} Given an element $\tau\in \pl$ which is different from the
identity, there exist two intervals $I$ and $J$ of the real line,
whose endpoints are dyadic integers, with $\tau(I)=J$, and such that
$I\cap J=\varnothing$.
\end{lemma}

\begin{proof}
The case when the slope of $\tau$ is always 1 or $-1$ is trivial. For a
map $t\mapsto t+k$ has a small interval (of length less than $k$)
whose image is disjoint from it. If $\tau=-Id$ the result is
trivial.

If the slope is not constantly equal to 1, it has a piece with slope
$\pm 2^i$ with $i\neq 0$. Assume without loss of generality (by possibly
taking $\tau^{-1}$ instead of $\tau$) that $i>0$.
Hence there are two intervals $[a,b]$ and $[c,d]$
such that $\tau(a)=c$ and $\tau(b)=d$ and also $d-c=2^i(b-a)$. It is
possible that $[a,b]$ and $[c,d]$ overlap, but since $[c,d]$ is much
larger than $[a,b]$ (at least twice the size), we can choose as $J$
a small interval inside $[c,d]$ which is disjoint from $[a,b]$. By
construction, the preimage $I$ of $J$ is in $[a,b]$, and hence $I$
and $J$ are disjoint.
\end{proof}

\begin{proof}[Proof of \fullref{thm:injective}]
 We now take a nontrivial $\tau\in\com{F}$. By the previous lemma, there exist intervals $I$ and $J$ satisfying the conditions stated
above and, in addition, that $I$, and hence $J$, have endpoints of the form 
$k/2^j$ and $(k+1)/2^j$. 
We consider all elements of $F$ whose support (that is, the
part where they are not the identity) is contained in $I$. Those
elements form a subgroup which is isomorphic to $F$ itself. Let $f$
be one such element. Since its support is inside $I$, its image
under the commensuration $\tau$, that is,
$f^\tau=\tau\circ f\circ \tau^{-1}$, has support inside $J$.

Hence, the distance (inside $F$) from $f$ to $f^{\tau}$ is given by
the distance from the identity to the element $f^\tau f^{-1}$. But
this element has its support inside the disjoint union $I \cup J$, and the two parts are
independent from each other (one given by $f$ and the other one by
$f^\tau$). By work of Cleary and Taback \cite{clearytaback}, 
this subgroup---elements with
support in $I \cup J$ which is a direct product of two clone subgroups in their terminology---is quasi-isometrically embedded in $F$. Hence,
we can take elements $f_n$ with support inside $I$ with arbitrarily
large norm, and hence $f_n^\tau f_n^{-1}$ has also arbitrarily large norm.
This proves that the image of $\tau$, a quasi-isometry, is not
at bounded distance from the identity and the proof is complete.
\end{proof}

%%%%%%%%%%%%%%%%%%%%   End of main body of article
%
%                             References
%
%   BiBTeX users uncomment the following line:
%

\bibliographystyle{gtart}

%

%\begin{thebibliography}
%\end{thebibliography}

\end{document}